\documentclass[12pt]{amsart}

\usepackage{amssymb,amscd,pb-diagram,mathrsfs}
\usepackage{fullpage}

\newcount\timehh\newcount\timemm
\timehh=\time 
\divide\timehh by 60 \timemm=\time
\def\Z{{   \mathbb Z }}
\def\Q{{\mathbb Q}}

\def\F{{\mathbb F}}
\def\1{{\bf 1}}

\def\u{ {\lambda}}
\def\v{ {\mu}}

\def\lcm{{\rm lcm}}

\def\ord{{\rm {ord}}}

\def\l{{\ell}}

\def\c2{\chi_2}

\def\hfl#1#2{\smash{\mathop{\hbox to
12mm{\rightarrowfill}}\limits_{\scriptstyle #1}^{\scriptstyle #2}}}

\def\Z{{   \mathbb Z }}
\def\Q{{\mathbb Q}}

\def\F{{\mathbb F}}

\def\l{\ell}

\def\La{L_n^{(\alpha)}(x)}

\newtheorem{theorem}{Theorem}[section]
\newtheorem{lemma}[theorem]{Lemma}

\newtheorem{definition}[theorem]{Definition}

\newenvironment{Proof}{\removelastskip\par\medskip
\noindent{\em Proof.} \rm}{\penalty-20\null\hfill$\square$\par\medbreak}

\title[Galois Group of Generalized Laguerre Polynomials]
{On the Galois Group of Generalized Laguerre Polynomials}

\author{
          Farshid Hajir\\
   {
    \protect \protect\sc\today\ --
    \ifnum\timehh<10 0\fi\number\timehh\,:\,\ifnum\timemm<10 0\fi\number\timemm
    \protect \, \, \protect 
  }
}
\email{hajir@math.umass.edu}
\address{Department of Mathematics \& Statistics, University of Massachusetts,
	Amherst, MA 01003-9318 USA}



\thanks{This work was supported by the National Science
Foundation under Grant No. 0226869.}
\begin{document}

\begin{abstract}
Using the theory of Newton Polygons, we formulate a simple criterion
for the Galois group of a polynomial to be ``large.''  For a fixed
$\alpha \in \Q - \Z_{<0}$, Filaseta and Lam have shown that the $n$th
degree Generalized Laguerre Polynomial $L_n^{(\alpha)}(x) =
\sum_{j=0}^n \binom{n+\alpha}{n-j}(-x)^j/j!$ is irreducible for all
large enough $n$.  We use our criterion to show that, under these
conditions, the Galois group of $\La$ is either the alternating or
symmetric group on $n$ letters, generalizing results of Schur for
$\alpha=0,1$.
%\vspace{.3in}
%
%\noindent {\sc Abstrait.} En utilisant la th\'eorie des polyg\^ones de
%Newton, on obtient un crit\`ere simple pour montrer que le groupe de
%Galois d'un polyn\^ome soit ``large.''  Si on fixe $\alpha \in \Q
%-\Z_{<0}$, Filaseta et Lam ont montr\'e que le Polyn\^ome Generalis\'e
%de Laguerre $L_n^{(\alpha)}(x) = \sum_{j=0}^n
%\binom{n+\alpha}{n-j}(-x)^j/j!$ est irreductible quand le degr\'e $n$
%est assez grand.  On utilise notre crit\`ere afin de montrer que,
%sous ces hypoth\`eses, le groupe de Galois de $\La$ est soit le groupe
%altern\'e, soit le groupe symmetrique, de degree $n$, g\'en\'eralisant
%des resultats de Schur pour $\alpha=0,1$.
\end{abstract}

\maketitle

%\begin{center}
%{\em \`A Georges Gras, \`a l'occasion de son 60\`eme anniversaire}
%\end{center}

\section{Introduction}

It is a basic problem of algebra to compute the Galois group of a
given irreducible polynomial over a field $K$.  If we order the monic
degree $n$ polynomials over $\Z$ by increasing height, then the
proportion which consists of irreducible polynomials with Galois group
$S_n$ tends to 1; for a more precise statement, see for example
Gallagher \cite{gallagher}.  Nevertheless, to prove that the Galois
group of a given polynomial is $S_n$ can be difficult if $n$ is large.
The algorithmic aspects of Galois group computations have witnessed a
number of recent advances, for which an excellent reference is the
special issue \cite{mmy} of the Journal of Symbolic Computation,
especially the foreward by Matzat, McKay, and Yokoyama.  Currently,
for rational polynomials of degree up to 15, efficient algorithms are
implemented, for instance, in {\sc gp-pari} and {\sc magma}.  An
important piece of any such algorithm is the collection of data
regarding individual elements of the Galois group, for which the
standard method is to factor the polynomial modulo various ``good''
primes (i.e. those not dividing its discriminant), obtaining the
cycle-type of the corresponding Frobenius conjugacy classes in the
Galois group.

Our first goal in this paper is to formulate a criterion which
exploits the properties of ``bad'' primes for proving that the Galois
group of a given polynomial is large.  The criterion is especially
efficacious if one suspects that a ``medium size'' prime (roughly
between $n/2$ and $n$) is wildly ramified in the splitting field of
the polynomial.  The criterion we give (Theorem \ref{coga}) follows
quite simply from the theory of $p$-adic Newton Polygons; it is used
in slightly less general form in Coleman \cite{coleman} and is
reminiscent of, but distinct from, a criterion of Schur \cite[\S
1]{sch67}.

Our second goal is to illustrate the utility of the criterion 
by using it to calculate the Galois group
for a certain family of polynomials, which we now introduce.  In the
second volume of their influential and classic work \cite{pz}, P\'olya
and Szeg\H o define the Generalized Laguerre Polynomial (GLP)
$$L_n^{(\alpha)}(x) = \sum_{j=0}^n \binom{n+\alpha}{n-j}
\frac{(-x)^j}{j!}.$$ The special case $\alpha=0$ had appeared much
earlier in the work of Abel \cite[p.~284]{abel} and Laguerre
\cite{laguerre}, and the general case can in fact be found in Sonin
\cite[p.~41]{sonin}.  Shortly after the publication of \cite{pz}, the
study of the algebraic properties of this family of orthogonal
polynomials was initiated by Schur \cite{sch67}, \cite{sch70}.

For instance, for the discriminant of the monic integral polynomial
$(-1)^n n! L_n^{(\alpha)}(x)$, we have the following formula of Schur
\cite{sch70}:
\begin{equation}
\label{disc} \Delta_n^{(\alpha)} =
\prod_{j=2}^n j^j (\alpha+j)^{j-1}.
\end{equation}
In particular, if $\alpha$ is not in $[-n,-2]\cap \Z$,
$L_n^{(\alpha)}(x)$ has no repeated roots.
%; we let $G_n^{(\alpha)}$ be its Galois group over $\Q$.  
For $\alpha=0,1$, Schur \cite{sch67}, \cite{sch70} established the
irreducibility of all $\La$ over $\Q$, and also showed that their
Galois groups are as large as possible, namely $A_n$ if
$\Delta_n^{(\alpha)}$ is a rational square, and $S_n$ otherwise.

%In \cite{sch67} and \cite{sch70}, Schur showed for $\alpha=0,1$ (classical
%Laguerre polynomial) and $\alpha=1$, that $L_n^{(\alpha)}(x)$ is
%irreducible in $\Q[x]$ for all $n$.  He also showed for these values
%of $\alpha$ that $G_n^{(\alpha)}$ always contains the alternating 
%group $A_n$.  By considering the value of $\Delta_n^{(\alpha)}$
%in $\Q^*/{\Q^*}^2$, he then computed: $G_n^{(0)}=S_n$ for all $n$,
%and $G_n^{(1)}=S_n$, unless $n$
%is odd or $n+1$ is an even square, in which case, $G_n^{(1)}=A_n$.

A number of recent articles on the algebraic properties of GLP have
appeared, including Feit \cite{feit}, Coleman \cite{coleman}, Gow
\cite{gow}, Filaseta-Williams \cite{fw}, Filaseta-Lam \cite{fl}, Sell
\cite{sell}, Hajir \cite{jnt}, \cite{hajir}, and Hajir-Wong \cite{hajir-wong}.
In particular, we have the following theorem of Filaseta and Lam
\cite{fl} on the irreducibility of GLP.  \par
\vspace{\baselineskip}
%\begin{theorem}[Filaseta-Lam]\label{fl}
\noindent {\bf Theorem.} (Filaseta-Lam) {\em
If $\alpha$ is a fixed rational number which is not a negative integer, 
then for all
but finitely many integers $n$,  
$L_n^{(\alpha)}(x)$ is irreducible
over $\Q$.}
%\end{theorem}
\par
\vspace{\baselineskip} In this paper, we provide a complement to
the theorem of Filaseta and Lam by computing the Galois group of $\La$
when $n$ is large with respect to $\alpha\in \Q-\Z_{<0}$.  Namely, we
prove the following result.
\begin{theorem}\label{mt}
Suppose $\alpha$ is a fixed rational number which is not a negative
integer.  Then for all but finitely many integers $n$, the Galois
group of $\La$ is $A_n$ if $\Delta_n^{(\alpha)}$ is a square and $S_n$
otherwise.
%\begin{enumerate}
%\item[i)]If the interval $((n+\alpha)/2,n-2)$ contains a prime, 
%and if $\La$ is irreducible over $\Q$, then its Galois group over $\Q$
%contains $A_n$.
\end{theorem}

\noindent {\bf Remarks.}  1. The hypothesis that $\alpha$ not be a
negative integer is necessary, as in that case, $\La$ is divisible by
$x$ for $n\geq |\alpha|$.  For a study of the algebraic properties of
$\La$ for $\alpha\in \Z_{<0}, n<|\alpha|$, see \cite{jnt}, \cite{sell}
and \cite{hajir}.

\noindent 2.  Using a different set of techniques, the following
companion to Theorem \ref{mt} is proved in \cite{hajir-wong}: If we
fix $n\geq 5$ and a number field $K$, then for all but finitely many
$\alpha\in K$, $\La$ is irreducible and has Galois group $A_n$ or
$S_n$ over $K$.  For each $n\leq 4$, infinitely many reducible
specializations exist, and for $n=4$, there are infinitely many
specializations which are irreducible but have $D_4$-Galois group,
cf. \cite[Section 6]{hajir}.

\noindent 3.  For integral $\alpha$, some cases where
$\Delta_n^{(\alpha)}$ is a square (giving Galois group $A_n$) are
\begin{enumerate}
\item[$\bullet$] $\alpha=1$ and $n\equiv 1\pmod{2}$ or $n+1$ is an odd square
(\cite{sch70}),
\item[$\bullet$] $ \alpha = n,$ and $n\equiv 2 \pmod{4}$ (\cite{gow}, it is not
yet known if all of these polynomials are irreducible \cite{fw}),
\item[$\bullet$] $ \alpha = -1-n$, and $n\equiv 0 \pmod{4}$
(\cite{sch67}, \cite{coleman}), \item[$\bullet$] $\alpha = -2-n$, and
$n\equiv 1 \pmod{4}$ (\cite{jnt}).  \end{enumerate} See \cite[\S
5]{hajir} as well as the above-cited papers for more details.

\noindent 4. The proofs of the Filaseta-Lam Theorem in
\cite{fl} and of Theorem \ref{mt} are both effective.

%is
%effective, i.e. one can specify a (large) bound $M(\alpha)$ so that
%$L_n^{(\alpha)}$ is irreducible over $\Q$ for $n\geq M(\alpha)$.  In
%this paper, we give an explicit function $N(\alpha)$ with the property
%that for $n\geq N(\alpha)$, the Galois group of $\La$ contains $A_n$
%assuming $\La$ is irreducible.  
%Filaseta and Lam do not give $M(\alpha)$
%explicitly, but it seems likely that $N(\alpha) \leq M(\alpha)$ since
%the former is linear, but the latter must almost certainly taken to 
%be exponential in $\alpha$ for the Filaseta-Lam proof to apply.

\section{A criterion for having large Galois group}

\subsection{Newton Polygons}

Let $K$ be a field equipped with a discrete valuation $v$ and a
corresponding completion $K_v$.  We assume $v$ is normalized,
i.e. $v(K^*)=\Z$, and employ the same letter $v$ to denote an extension
of this valuation to an algbraic closure $\overline{K_v}$ of $K_v$.

For a polynomial $f(x)=a_n x^n + a_{n-1}x^{n-1} + \cdots + a_0\in
K[x]$ with $a_0 a_n \neq 0$, the $v$-adic Newton Polygon of $f(x)$,
denoted $NP_v(f(x))$, is defined to be the lower convex hull of the
set of points
$$S_v(f) = \{(0,v(a_0)), (1, v(a_1)), \cdots, (n,v(a_n))\}.$$  It is the
highest polygonal line passing on or below the points in $S_v(f)$.
The points where the slope of the Newton polygon changes
(including the rightmost and leftmost points) are called the
{\em corners} of $NP_v(f)$; their $x$-coordinates are the {\em breaks}
of $NP_v(f)$.

For the convenience of the reader, we recall the main theorem about
$v$-adic Newton Polygons as well as a useful corollary due to Coleman
\cite{coleman}.  See, for instance, Gouv\^ea \cite{gouvea}.
%(Newton had introduced the concept for the
%field of Laurent series).  
A very nice survey of the uses of the
Newton Polygon for proving irreducibility is Mott \cite{mott}.  For 
generalizations to several variables, see Gao \cite{gao} and references
therein.

\begin{theorem}[Main Theorem of Newton Polygons]
\label{NP} Let $(x_0,y_0),(x_1,y_1),\ldots,(x_r,y_r)$ denote the
successive vertices of $NP_v(f(x)).$ Then there exist polynomials
$f_1,\ldots,f_r$ in $\Q_p[x]$ such that\begin{enumerate}\item[i)]
$f(x)=f_1(x)f_2(x)\cdots f_r(x),$\item[ii)]the degree of $f_i$ is
$x_i-x_{i-1},$\item[iii)]all the roots of $f_i$ in
$\overline{K_v}$ have $v$-adic valuation
$-(y_i-y_{i-1})/(x_i-x_{i-1}).$\end{enumerate}
\end{theorem}

\subsection{Newton Index}

We now suppose that $K$ is a global field, i.e. $K$ is a finite
extension of $\Q$ (number field case) or of $\F(T)$, where $\F$ is a
finite field (function field case).  A global field $K$ enjoys the
property that for a given element $\alpha \in K$, $v(\alpha)=0$ for
all but finitely many valuations $v$ of $K$.

\begin{definition}
Given $f\in K[x]$, the {\em Newton Index} of $f$ ${\mathscr{N}}_f$ 
is defined to be the least common multiple of the
denominators (in lowest terms) of all slopes of $NP_v(f)$
as $v$ ranges over all normalized discrete valuations of $K$.
\end{definition}
Note that 0 is defined to have denominator 1, so slope $0$ segments of
$NP_v(f)$ do not contribute to ${\mathscr N}_f$.  On the other hand,
for all but finitely many $v$, the coefficients of $f$ all have
$v$-adic valuation $0$ so $NP_v(f)$ consists of a single slope 0
segment.  Hence, ${\mathscr {N}}_f$ is well-defined and effectively
computable (for monic $v$-integral polynomials, we need only compute
the Newton Polygon for those valuations that do not vanish on the
constant coefficient).  It is clear that $\mathscr{N}_f$ is a divisor
of $\lcm(1,2,\ldots,n)|n!$ where $n$ is the degree of $f$; the latter
property in fact holds for an arbitrary field $K$, so the Newton index
is well-defined for any field $K$, though possibly not in the sense
that it is necessarily effectively computable.

We now formulate a criterion for an irreducible polynomial to have
``large'' Galois group.  
%The criterion is useful when the ramification in $\Q[x]/(f)$ is ``generic.''  
The key idea appears in Coleman's
computation \cite{coleman} of the Galois group of the $n$th Taylor
polynomial of the exponential function, which incidentally is the GLP
$(-1)^n L_n^{(-1-n)}(x)$.

\begin{theorem}\label{coga}
Suppose $K$ is a global field and $f(x)$ is an irreducible polynomial
in $K[x]$. Suppose $g(x)=f(x-\mu)$ for some $\mu \in K$.  Then
${\mathscr N}_g$ divides the order of the Galois group of $f$ over $K$.
Moreover, if ${\mathscr N}_g$ has a prime divisor $q$ in the range
$n/2 < q < n-2$, where $n$ is the degree of $f$, then the Galois group
of $f$ contains $A_n$.
\end{theorem}
\begin{Proof}
Suppose $v$ is a valuation of $K$ and $q$ is an arbitrary divisor of
the denominator of some slope $s$ of the $v$-adic Newton polygon of
$g$.  Clearly, $f$ and $g$ have the same splitting field and the same
Galois group.  It suffices to show that $q$ divides the order of the
Galois group of $g$ over $\Q$.  By the main theorem of Newton polygons
\ref{NP}, there exists a root $\alpha \in \overline{K}_v$ of $g$ with
valuation $-s$.  Since $q$ divides the denominator of $s$, $q$ divides
the ramification index $e$ of $K_v(\alpha)/K_v$.  But $e$ divides the
degree $[K_v(\alpha):K_v]$, which in turn divides the order of the
Galois group of $g$ over $K_v$, hence also over $K$.  If $q$ is a
prime in the interval $(n/2,n-2)$, then the Galois group of $g$
contains a $q$-cycle, so it must contain $A_n$ by a theorem of Jordan
\cite{jordan}
(or see, for instance, Hall's book \cite[Thm~5.6.2 and 5.7.2]{hall}).
\end{Proof}

\noindent {\bf Remark.}  Schur proved a similar result (\cite[\S 1,
III]{sch67}), namely, if the discriminant of a number field $K$ of
degree $n$ is divisible by $p^n$, then the Galois closure $L$ of $K$
has degree $[L:\Q]$ divisible by $p$.  In general, if $p$ divides the
discriminant of an irreducible polynomial $f$, it is not easy to
determine the $p$-valuation of the discriminant of the stem field
$\Q[x]/(f)$; thus, each of Theorem \ref{coga} and Schur's criterion
can be useful depending on whether we have information about the
discriminant of the field or that of the defining polynomial.  Neither
criterion is useful when the discriminant of $f$ is square-free, for
example, since in that case, all the non-trivial ramification indices
are $2$.  On the other hand, over base field $\Q$, irreducible
polynomials with square-free discriminant also have Galois group $S_n$
see e.g. Kondo \cite{kondo}; the proof of this fact uses the
triviality of the fundamental group of $\Q$.

\section{Proof of Theorem \ref{mt}}

We now let $K=\Q$.  For a prime $p$, we write $NP_p$ in place of
$NP_v$ where $v=\ord_p$ is the $p$-adic valuation of $\Q$.

\begin{lemma}\label{key}
Let $f(x)= \sum_{j=0}^n \binom{n}{j} c_j x^j \in \Q[x]$ be an
irreducible polynomial of degree $n$ over $\Q$.  Suppose there exists a 
prime $p$ satisfying
\begin{enumerate}
\item[i)] $n/2 < p < n-2$,
\item[ii)] $\ord_p(c_j)\geq 0$ for $0\leq j \leq n$,
\item[iii)] $\ord_p(c_j)=1$ for $1\leq j \leq n-p$,
%is positive but not divisible by $p$,
%\item[iv)] $\ord_p(c_j) \geq \ord_p(c_0)$ for $1\leq j \leq n-p$
\item[iv)] $\ord_p(c_p)=0$.
\end{enumerate}
Then the Galois group of $f$ over $\Q$ contains $A_n$.  
\end{lemma}
\begin{Proof}
It is easy to check that $\binom{n}{j}$ is divisible by $p$ if and
only if $n-p+1 \leq j \leq p-1$.  The given assumptions then guarantee
that $(0,1)$ and $(p,0)$ are the first two corners of $NP_p(f)$.
Therefore, $-1/p$ is a slope of $NP_p(f)$, hence $p|{\mathscr N}_f$
and we are done by Theorem
\ref{coga}.
\end{Proof}
%Since Filaseta and Lam \cite{fl} have shown that $\La$ is irreducible
%over $\Q$ for $n$ large with respect to $\alpha$ (assuming $\alpha$ is
%not a negative integer), the proof of Theorem \ref{mt} will be
%complete once we establish the following theorem.
%\begin{theorem}\label{compute}
%Suppose $\alpha\in \Q$ is not a negative integer.  We write
%$\alpha=\u/\v$ in lowest terms, i.e. with $\gcd(\u,\v)=1$ and $\v\geq
%1$.  There exists an effectively computable constant $M(\u,\v)$ such
%that for $n\geq M(\u,\v)$, $\La$ is irreducible over $\Q$, with Galois
%group either $A_n$ or $S_n$.
%\end{theorem}

We are now ready to prove the Main Theorem.

\begin{proof}[Proof of Theorem \ref{mt}]
We write
$\alpha=\u/\v$ in lowest terms, i.e. with $\gcd(\u,\v)=1$ and $\v\geq
1$.  By assumption, $\alpha$ is not a negative integer.
We will work with the normalized (monic, integral) polynomial
$$f(x):=\v^n n! L_n^{(\u/\v)}(-x/\v) = \sum_{j=0}^n
\binom{n}{j} (n\v+\u)((n-1)\v+\u)\cdots ((j+1)\v+\u) x^j.$$ 
We wish to apply
Lemma \ref{key} to it, so we let \begin{equation}\label{cj}
c_j = \prod_{k=j+1}^n (k\v + \u), \qquad 0\leq j \leq n,
\end{equation}
and seek a suitable prime $p$ satisfying the conditions of the Lemma.

By a suitably strong form of Dirichlet's theorem on primes in
arithmetic progressions, there exists an effective constant $D(\v)$
such that if $x\geq D(\v)$ and $h\geq x/(2\log^2 x)$, the interval
$[x-h,x]$ contains a prime in the congruence class $\u \bmod{\v}$ (see
Filaseta-Lam \cite[p.~179]{fl}).
Taking $x=n-3\geq D(\mu)$, we 
find that 
for some integer $\l \in [1,n]$, $p=\v\l+\u$ is a prime satisfying
\begin{equation}\label{p}
\frac{n\v+\v+\u}{\v+1} \leq p \leq n-3,
\end{equation}
as long as  $$
\frac{1-3/n}{2\log^2(n-3)} + \frac{3 + \lambda/(\v+1)}{n} \leq \frac{1}{\v+1},
$$
which clearly holds for all $n$ large enough with respect to $\u,\v$.

We now fix a prime $p=\v\l+\u$ satisfying (\ref{p}).
For such a prime $p$, let us check the hypotheses of Lemma \ref{key}.
%\begin{enumerate} 
%\item[i)] 
We have $(n\v+\v+\u)/(\v+1)>n/2$ if and only
if \begin{equation}\label{te} 
n(\v-1)>-2\v-2\u.
\end{equation}  
Since $\alpha$ is not a negative integer,
if $\v=1$, then $\lambda\geq 0$, so (\ref{te}) holds for all $n$.
If $\v>1$, we simply need to take $n> -2(\v+\u)/(\v-1)$ in order to
achieve
$n/2<p<n-2$, giving us i).
 Our $c_j$ are integral so ii) holds trivially.  Before we 
discuss iii),
let us note that in the congruence class
$\u \bmod \v$,  the smallest multiple of $p$
larger than $p$ is $(\v+1)p$, and,
similarly, the largest multiple of $p$ in this congruence class which
is less than $p$ is $(-\v+1)p$.  Now we claim that, for $n$ large enough,
we have
\begin{eqnarray}\label{b1}
(-\v+1)p &<& \u+\v
\end{eqnarray}
as well as \begin{eqnarray}\label{b2}
\u + \v n &<& (\v+1)p. 
\end{eqnarray}
Indeed, if $\mu=1$, then (\ref{b1}) holds for all $n$, while for
$\mu\geq 2$, $n\geq -2\u$ implies (\ref{b1}); 
%moreover, (\ref{p})
%implies $$ p> \frac{n\v+\u}{\v+1},$$ which is easily seen to be
%equivalent to (\ref{b2}).  
moreover, (\ref{b2}) is a direct consequence of (\ref{p}).
We conclude from (\ref{b1}) and (\ref{b2}) that
$\ord_p(c_0)=1$.  From (\ref{cj}), we then read off that $\ord_p(c_j)=
1$ for $0\leq j \leq \l -1$, and $\ord_p(c_j)=0$ for $\l\leq j\leq n$.
One easily checks that (\ref{b1}) and (\ref{b2}) give exactly 
$p > \l -1$ and $n-p<\l$, i.e. iii) and iv).  By Filaseta-Lam \cite{fl},
there is an effectively computable constant $N(\alpha)$ such
that $f(x)$ is irreducible for $n\geq N(\alpha)$.
Thus, all the conditions of 
Lemma \ref{key} hold, and the proof of the theorem is complete.
\end{proof}
\noindent{\bf Remark.} Note that the proof simplifies in the case 
where $\alpha$ is a non-negative integer, giving: If 
$L_n^{(\alpha)}(x)$ is irreducible and if
there is a prime
$p$ in the interval $( (n+\alpha)/2,n-2)$, then the Galois group of 
$L_n^{(\alpha)}(x)$ contains $A_n$.  By 
\cite[Corollary 3.2]{hajir}, 
the specified interval contains
a prime as long as
$n\geq \max(48-\alpha,8+5\alpha/3)$.

\end{document}